\documentclass[journal]{IEEEtran}

\usepackage{amsmath,amsfonts,amsthm}
\usepackage{algorithm,algorithmic}
\usepackage{subcaption}
\usepackage{textcomp}
\usepackage{stfloats}
\usepackage{verbatim}
\usepackage{graphicx}
\usepackage{epstopdf}
\usepackage{balance}
\usepackage{placeins}

\usepackage[hidelinks]{hyperref}

\graphicspath{{figures/}}
\newtheorem{theorem}{Theorem}
\newtheorem{lemma}{Lemma}

\newtheorem{remark}{Remark}

\begin{document}

\title{Distributed Attacker Assignment and Interception Control in Multi-Agent Target-Attacker-Defender Games}

\author{Pengyang Fan, Shaoqing Zhang, Chao Zhai, and Weiguo Xia%
\thanks{Pengyang Fan and Weiguo Xia are with the Key Laboratory of Intelligent Control and Optimization for Industrial Equipment of Ministry of Education and School of Control Science and Engineering, Dalian University of Technology, Dalian 116024, China (e-mail: \href{mailto:fanpengyang@mail.dlut.edu.cn}{fanpengyang@mail.dlut.edu.cn}; \href{mailto:wgxiaseu@dlut.edu.cn}{wgxiaseu@dlut.edu.cn}). Chao Zhai is with the School of Artificial Intelligence and Automation, Hubei Key Laboratory of Advanced Control and Intelligent Automation for Complex Systems, and Engineering Research Center of Intelligent Technology for Geo-exploration, Ministry of Education, China University of Geosciences, Wuhan 430074, China (e-mail: \href{mailto:zhaichao@amss.ac.cn}{zhaichao@amss.ac.cn}). Shaoqing Zhang is with the Key Laboratory of Intelligent Control and Optimization for Industrial Equipment of Ministry of Education and School of Control Science and Engineering, Dalian University of Technology, Dalian 116024, and also with Shenyang Aircraft Design and Research Institute, Shenyang 110035, China (e-mail: \href{mailto:shaoqing_c_z@126.com}{shaoqing\_c\_z@126.com}).}}

\maketitle

\begin{abstract}
This paper addresses the Target-Attacker-Defender (TAD) game within a planar environment, focusing on the protection of a stationary target against incoming threats. To compensate for potential speed disadvantages of the defenders, we leverage geometric pre-positioning within an annular defense perimeter. The problem is formulated as an area interception problem, for which we propose a distributed multi-agent algorithm consisting of two primary modules: target allocation and interception control. For target allocation, an equitable region partitioning strategy is employed to divide the defense area into annular sectors and assign defenders to attackers. Under specified angular-separation and parameter conditions, the proposed partitioning mechanism ensures that each defender is assigned at most one attacker, thereby reducing the complexity of multi-agent assignment. For the control layer, movement protocols based on zero-sum differential games and coverage optimization are designed, achieving Nash equilibrium strategies while maintaining effective barrier coverage. Numerical simulations show that the proposed method enables slower defenders to intercept faster attackers through proactive deployment and maintains nearly constant local computation time as the swarm size increases. 
\end{abstract}


	\section{Introduction}

The protection of critical infrastructures, such as chemical plants~\cite{argenti2017study} and airports~\cite{lykou2020defending}, has driven extensive research into the coordination of multi-robot teams for area security. A standard framework for this task is multi-agent barrier coverage~\cite{luo2016strong,li2011study,SHI201595}, which focuses on deploying and coordinating agents along target boundaries for perimeter monitoring. Nevertheless, traditional coverage-control strategies often lack active interception capabilities or rely on predefined patrolling patterns~\cite{cheng2009distributed,cheng2011decentralized}. While stochastic approaches, such as Kalman filtering~\cite{zhai2016coverage}, can estimate intruder distributions, they remain vulnerable to intelligent attackers who exploit defensive motion patterns to design evasive trajectories. 

To counter intelligent attackers, recent research has increasingly turned to game-theoretic formulations, specifically Target-Attacker-Defender (TAD) games and Reach-Avoid games. TAD games are differential games in which attackers seek to capture a target, whereas defenders attempt to intercept them. Reach-Avoid games, by contrast, model defenders and obstacles as the avoidance set and the target region as the reach set~\cite{margellos2011hamilton}. Both TAD games and Reach-Avoid games can be traced back to Isaacs' seminal work on pursuit-evasion problems~\cite{isaacs1999differential}. Although pursuit-evasion differential games have been widely studied~\cite{weintraub2020introduction,patsko2011homicidal,li2008stochastic}, the presence of a target in TAD games introduces additional analytical challenges~\cite{rusnak2011guidance,garcia2020defense,liang2020differential,liang2022reconnaissance}.

Early studies in this area primarily considered single-attacker-single-defender scenarios, and more recent work has shifted toward multi-attacker-multi-defender settings. Although interception strategies involving multiple defenders against a single attacker have been investigated~\cite{singh2022dynamic}, these methods become less effective as the number of attackers increases. Consequently, target assignment emerges as a pivotal challenge in multi-agent scenarios, requiring	 optimal pairings between defenders and attackers. 

Early approaches to this problem mainly rely on assignment-cost-based allocation methods. Chen et al.~\cite{chen2014multiplayer} adopted a matching-based approach, which was further extended in later studies~\cite{chen2017pairwise}. More recently, Manoharan et al.~\cite{manoharan2023multi} combined neural networks with model predictive control; however, such methods typically depend on a centralized controller or a leader agent for assignment and therefore lack a distributed framework. Soni et al.~\cite{Soni2025dynamic} proposed a dynamic assignment strategy within a Linear Quadratic Differential Game (LQDG) framework, in which the Euclidean distance between defenders and attackers is minimized as a dynamic cost. However, accurately evaluating assignment costs in complex scenarios in real time remains computationally challenging. Moreover, in large-scale attack scenarios, allocation methods based solely on interception time may be inadequate, because they neglect the critical impact of the global spatial distribution of agents. As a result, defenders may become excessively concentrated in local areas, thereby undermining the integrity of the overall defensive barrier.

Another major class of target allocation approaches relies on game-theoretic reachability analysis. Yan et al.~\cite{yan2020task} derived explicit geometric boundaries to determine interception feasibility, thereby transforming the assignment problem into an integer programming problem. Based on this, matching-based capture strategies for three-dimensional heterogeneous games are proposed in~\cite{yan2022matching}, where a maximum matching algorithm is used to optimize assignments within a reach-avoid framework. Coon et al.~\cite{coon2017control} employed isochrone-based analysis to assess interception feasibility and used interception time as the cost metric to construct a bipartite graph for matching. Similarly, Li et al.~\cite{Li2025Distributed} combined isochrone analysis with winning-region characterization to address heterogeneous games involving players with damped double-integrator dynamics. They developed a distributed allocation algorithm that optimizes the pursuer team's makespan through an iterative procedure combining bipartite matching and binary search. However, such geometric approaches generally rely on the restrictive assumption that defenders are faster than attackers, which limits their applicability against high-speed adversaries. A related line of work adopts incentive-driven implicit allocation, where task priorities are embedded into the payoff structure so that agents associated with higher weights or rewards are more likely to be selected through equilibrium computation. For example, Wang et al.~\cite{wang2025multi} addressed the problem of limited perception in multi-agent TAD games through an inverse game-theoretic framework. However, this type of implicit allocation mechanism depends heavily on tuning incentive parameters and weighting coefficients, which may make parameter calibration difficult in large-scale scenarios.

In summary, existing target assignment and interception approaches for active defense face two critical bottlenecks. First, most current strategies assume that defenders have a speed advantage. In many practical engagements, however, defenders may be slower than attackers, rendering direct pursuit ineffective. We argue that defenders can compensate for their speed disadvantages by leveraging strategic pre-positioning to establish geometric advantages and thereby achieve successful interception within a barrier-defense framework. Second, in large-scale attack scenarios, traditional pairing and assignment algorithms often incur prohibitive computational overhead and require sensitive parameter tuning, which limits their responsiveness in dynamic environments.

To address these challenges, this paper proposes a distributed barrier interception framework for multi-agent TAD games. The proposed approach integrates the geometric advantages of barrier coverage with game-theoretic decision-making to improve defensive performance against faster adversaries while maintaining computational efficiency in large-scale deployments. The main contributions of this paper are summarized as follows:

\begin{enumerate}
    \item A distributed target allocation algorithm based on equitable region partitioning is developed to construct attacker-defender pairs, ensuring that each defender is assigned at most one attacker under specified conditions.
    \item Coverage performance functions and game cost functions are introduced to guide both idle and engaged defenders, thereby ensuring the coverage quality and defensive effectiveness of the barrier interception process.
    \item A distributed barrier interception algorithm is proposed, which provides a theoretical foundation for the allocation mechanism and derives optimal strategies for defenders through zero-sum differential game analysis.
\end{enumerate}

The remainder of this paper is organized as follows. Section~\ref{sec:4.1} presents the problem formulation, Section~\ref{sec:4.2} provides the theoretical analysis, Section~\ref{sec:4.3} reports the numerical simulation results, and Section~\ref{sec:4.4} concludes the paper.
	
	\section{Problem Modeling}\label{sec:4.1}
	This section models the multi-agent barrier interception problem. First, mathematical modeling and the game model of barrier interception are introduced, followed by the proposition of termination conditions for the game. Finally, optimization formulations for defenders and attackers are presented.
	
	\subsection{Modeling of Barrier Interception}
	
\begin{figure}
		\centering
		\includegraphics[width=0.9\linewidth]{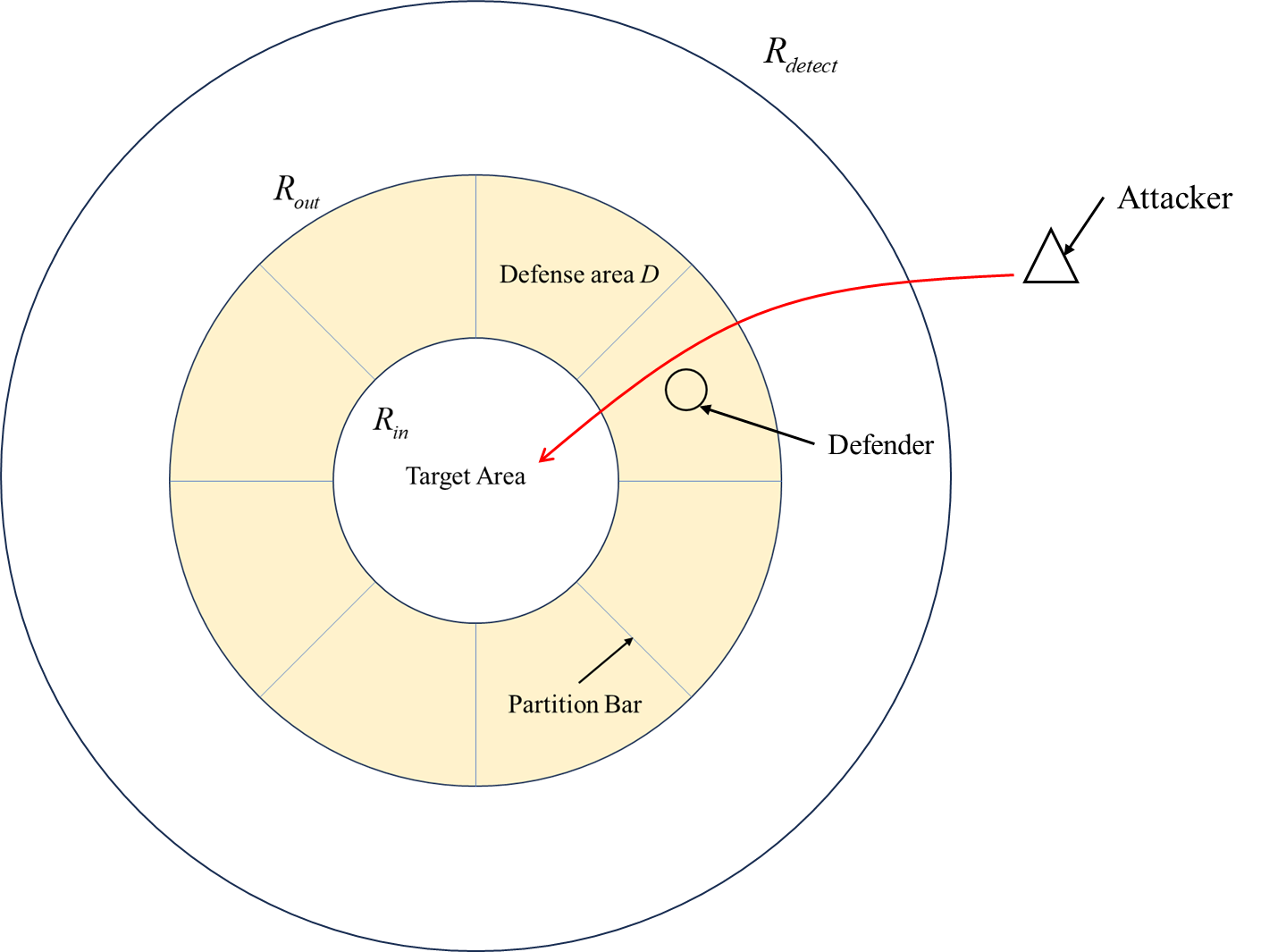}
		\caption{Diagram illustrating barrier interception: Defender (represented by circle) positioned on the circular barrier, Attacker (represented by triangle) outside the circular barrier, with partition points (depicted by black dots) on the barrier. The attacker aims to reach the target located at the center while the defender attempts to intercept the attacker's progress.}
		\label{fig:shiyitu01}
	\end{figure}

	Consider the defense area $D$, which is an annulus centered at the origin $O$, where the target is located. The defense area is bounded by an inner protection radius $R_{in}$ and an outer engagement radius $R_{out}$. Additionally, we define a detection radius $R_{detect}$, where $R_{detect} > R_{out}$. The region between $R_{out}$ and $R_{detect}$ serves as an early warning zone. As shown in Fig.~\ref{fig:shiyitu01}, defenders are positioned on the circular barrier between radius $R_{in}$ and $R_{out}$, while attackers are located outside $R_{detect}$, aiming to reach the target at the center. Defenders aim to intercept the attackers before they reach the target.
	
	Consider the scenario with $N$ defenders and $M$ attackers, where the capture radius of defenders is $R_d^c\in\mathbb{R}^+$. We consider the situation where the number of defenders exceeds that of attackers, i.e., $N>M$.
	Defenders are represented by the set $I_d=\{1,2,\ldots,N\}$, and attackers by $I_M=\{1,2,\ldots,M\}$.
	The defense area $D$ is partitioned into $N$ sub-regions, $E_1, \ldots, E_N$, by $N$ radial dividers. Each defender $i$ is responsible for the sub-region $E_i$.
	Moreover, the Cartesian coordinate vector for defenders is denoted by $p_d = [p_d^1 , p_d^2 , \cdots,  p_d^N]\in \mathbb{R}^{2\times N}$, with $p_d^i\in E_i$ for $i\in I_d$. The Cartesian coordinates for attackers are denoted by $p_a = [p_a^1 , p_a^2 , \cdots,  p_a^M]\in \mathbb{R}^{2\times M}$, with $p_a^j\in \mathbb{R}^2$ for $j\in I_M$.
	
	{We assume that the defenders possess a limited sensing range. An attacker $j$ is detected and its phase information becomes available to the defender network only when it enters the detection zone, i.e., $R_{in} \leq \|p_a^j\| \leq R_{detect}$. However, the active interception (engagement) is only permitted when the attacker enters the engagement zone, i.e., $\|p_a^j\| \leq R_{out}$. Attackers outside $R_{detect}$ remain undetectable.}
	
	The phase angle vector for defenders is denoted by $\varphi_d$, where $\varphi_d = (\varphi_d^1,\varphi_d^2,\ldots,\varphi_d^N)\in\mathbb{R}^N$. The numbering of phase angles follows the rule \(0\leq\varphi_{d}^1(0)<\varphi_{d}^2(0)<\cdots<\varphi_{d}^N(0)<2\pi \). The communication topology for defenders can be represented by an undirected cycle graph $\mathcal{G}_d = (I_d,\mathcal{E}_d)$, where $\mathcal{E}_d = \{(i, i+1)\mid i = 1, 2, \ldots, N-1\}\cup\{(N,1)\}$.
	The phase angle for attackers is denoted by $\varphi_a$, where $\varphi_a = (\varphi_a^1,\varphi_a^2,\ldots,\varphi_a^M)\in\mathbb{R}^M$.	
	
	We utilize the region partition method to divide the annular defense area into several subregions. We divide the annulus using $N$ virtual partition rays radiating from the origin. These rays are defined by their angles $S$, where $S = \{s_1, s_2,\ldots,s_N\}$ and $s_i \in [0, 2\pi)$ represents the angle of the partition ray for the $i$-th defender.
	The region covered by the $i$-th defender is denoted by $E_i$, which corresponds to the annular sector bounded by the partition angles $s_i$ and $s_{i+1}$.
	In addition, the partition angle $s_{i}$ should be adjusted to a suitable value to allocate attackers reasonably within the angular interval $[s_{i}, s_{i+1})$.
	
	Let $\rho(\theta)$ represent the spatial distribution of attackers invading at phase $\theta$, which is fitted by the defender network.
	The workload associated with region $E_i$ and defender $i$ is
	defined as
	\begin{equation}
		m_i = \int_{s_i}^{s_{i+1}} \rho(\theta) d\theta.
		\label{eq:bi_m}
	\end{equation}
	
	The following subsection will discuss how to model the spatial distribution of attackers and determine the placement of partition rays.
	
	\subsection{Density Function Fitting and Region Partition Dynamics}
	Prior research suggests that under single-integrator dynamics, the optimal strategy for attackers is to advance directly towards the target. When defenders lie on the path of attackers, attackers move along the tangent to the capture radius of defenders until the tangent aligns with the direction towards the target. At this point, attackers move away along the tangent, thereby avoiding capture and attacking the target~\cite{venkatesan2014target}. Although this pertains to single-integrator dynamics, it provides insights into our research, implying that attackers essentially do not alter their phases since the target is located at the origin. Hence, employing Gaussian kernel density estimation to fit the spatial distribution of attackers based on their phase angles is reasonable.
	
	Since the target is at the origin, the phase angle is the critical feature for interception. Note that this phase information is unknown before the attackers invade the detection area. Therefore, the density function fitting relies solely on the real-time phase data of attackers currently located within the detection range.
	
	For the defender network, let $\mathcal{I}_{obs}(t) = \{j \in I_M \mid R_{in} \le \|p_a^j(t)\| \le R_{detect}\}$ denote the set of observable attackers at time $t$. The phase angles of these observable attackers serve as sample data in Gaussian density estimation. We define the dynamic density function as follows:
	\begin{equation}
		\rho(\theta, t) = \sum_{j\in \mathcal{I}_{obs}(t)} K(\theta,\varphi_a^j,h),
		\label{eq:rho_e}
	\end{equation}
	where $K$ is the kernel, a non-negative function, and $h>0$ is the scaling factor. In contrast to traditional kernel density estimation, Equation~(\ref{eq:rho_e}) is not normalized by the total number of samples, which better determines the estimated number of attackers in the environment.
	
	Due to the periodicity of phase angles, the kernel function is redesigned based on the Gaussian kernel as follows:
	\begin{equation}
		K(\theta,\varphi,h)= e^{-h {\delta^2( \varphi,\theta)}},
		\label{eq:ke_gus}
	\end{equation}
	where $\delta(\varphi,\theta):[0,2\pi)\times[0,2\pi)\longrightarrow[-\pi,\pi)$ denotes the signed shortest angular difference from $\theta$ to $\varphi$, defined as positive in the clockwise direction and negative in the counterclockwise direction. It is evident that Equation~(\ref{eq:ke_gus}) is an even function, thus yielding the following inequality:
	\begin{equation}
		\int_0^{2\pi} K(\theta,\varphi,h) d\theta = 2\int_0^\pi e^{-hx^2}dx < \sqrt{\frac{\pi}{h}}.
		\label{eq:hpi}
	\end{equation}
	Similarly, an inequality related to the total number of attackers can be derived:
	\begin{equation*}
		\int_0^{2\pi} \rho(\theta) d\theta =  2M\int_0^\pi e^{-hx^2}dx < M \sqrt{\frac{\pi}{h}}.
	\end{equation*}
	This implies that the defender network fitted with the improved Gaussian kernel can more accurately determine the quantity and spatial distribution of the attacker network through the scaling factor $h$. For convenience, let
	$$
	\mathbf{e}(x)=\int_0^x e^{-h\eta^2}d\eta.
	$$
	It is observed that $\mathbf{e}: \mathbb{R} \to (-\sqrt{\frac{\pi}{4h}}, \sqrt{\frac{\pi}{4h}})$.
	Furthermore, $\lim_{{x \to +\infty}} \mathbf{e}(x) = \sqrt{\frac{\pi}{4h}}, \quad \lim_{{x \to -\infty}} \mathbf{e}(x) = -\sqrt{\frac{\pi}{4h}}$.

	Referring to the region partitioning algorithm designed in~\cite{zhai2019distributed}, to ensure that the defender network obtains a balanced workload, the following dynamic equation for partitioning angle $s_i$ is designed:
	\begin{equation}
		\dot s_i = \kappa_s (m_{i}-m_{i-1}) ,\quad i\in I_d,
		\label{eq:s}
	\end{equation}
	where $\kappa_s$ is a positive constant, and $m_0=m_N$. When the evolution of $s_i$ exceeds the interval $[0,2\pi)$, the modulo operation is applied resulting in $s_i = \mod(s_i,2\pi)$, which returns the remainder of $s_i$ divided by $2\pi$.
	It should be noted that to achieve dynamic partitioning in~(\ref{eq:s}), an agent needs to periodically collect information about the attacker network and broadcast the fitted density function to the defender network. Due to the requirements of interception applications, this period should be relatively short.
	
	However, relying on a centralized node to broadcast information across a large defense area may create a single point of failure. To overcome this drawback, we further propose a fully distributed density-fitting mechanism. Let \(I_M^i\) denote the set of attackers located within the annular sector \(E_i\), namely, the attackers monitored by the \(i\)-th defender.

	For the $i$-th defender, its Gaussian-fitted density function is defined as follows:
	\begin{equation}
		\rho_i^*(\theta) = \sum_{j\in I_M^i} K(\theta,\varphi_a^j,h).
		\label{eq:rho8}
	\end{equation}
	Agents need to be aware of the attacker information within their subregions and require communication to obtain information about the attacker network. Thus, Algorithm~\ref{tab:algong} is designed to enable defenders to extract relevant information about attackers by leveraging the communication topology $\mathcal{G}_d$.
	\begin{algorithm}
		\caption{\label{tab:algong}Attack Information Fusion Algorithm}
		\hspace*{0.02in} For any $i\in I_d$, the $i-$th defender agent performs as follows:
		\begin{algorithmic}[1]
			\STATE Detect attackers in local annular sector $E_i$
			\STATE Compute local attacker density function $\rho_i^*(\theta)$
			\STATE Save attacker information: $G_i \gets \{\rho_i^*(\theta)\}$
			\STATE Receive $G_j$ from neighbors based on $\mathcal{G}_d$
			\WHILE{$G_i \neq G_i \bigcup G_j$}
			\STATE Broadcast $G_i$ to neighbors based on $\mathcal{G}_d$
			\STATE Receive $G_j$ from neighbors based on $\mathcal{G}_d$
			\STATE Merge local and received attacker information: $G_i \gets G_i \cup G_j$
			\ENDWHILE
			\RETURN Merged attacker information set $G_i$ for agent $i$
		\end{algorithmic}
	\end{algorithm}
	
	Algorithm~\ref{tab:algong} models local attacker information using local detection and Gaussian-fitted density functions. Through the communication topology $\mathcal{G}_d$, agents can propagate and integrate local information. Through iterations, agents gradually integrate information from neighbors to achieve the fusion of global attacker information. Let $G_i$ represent the consolidated set of attacker density components acquired by defender $i$ upon the convergence of Algorithm~\ref{tab:algong}. Consequently, the local Gaussian density function can be formulated as follows:
	\begin{equation}
		\rho_i(\theta) = \sum_{\rho^* \in G_i}  \rho^*(\theta).
		\label{eq:rhoi}
	\end{equation}
	Ideally, the local density function $\rho_i(\theta)$ for defender $i$ converges to the global density function $\rho(\theta)$.
	
	Through the proposed density-fitting and region-partitioning algorithms, the originally complex multi-agent interception problem can be effectively decoupled. As will be rigorously shown in Section~\ref{sec:4.2}, this dynamic partitioning mechanism ensures that each defender's subregion contains at most one assigned attacker. Consequently, the macroscopic multi-agent engagement can be reduced, in a mathematically tractable manner, to a collection of localized one-to-one barrier interception games. On this basis, the game model presented in the following subsection concentrates on the interaction between a single attacker and a single defender.

	\subsection{Agent model and differential game}
	Consider an attacker $j$ and a defender $i$. The attacker $j$ can either enter the target region from the defender's assigned region $E_i$ or move to another defender's region with no attackers present (as in Theorem~\ref{se4:th2}).
	The defender is capable of moving freely within its assigned sub-region $E_i$ to intercept the attacker, while the attacker can move freely across the entire 2D plane.
	Their control inputs are defined based on single-integrator kinematics as follows:
	\begin{equation}\label{eq:ud}
		\dot{p}_d^i = u_d^i, \quad
		\dot{p}_a^j = u_a^j,
	\end{equation}
	where $u_d^i \in \mathbb{R}^2$ and $u_a^j \in \mathbb{R}^2$. The state of the game system is denoted by $x$, where $x := ({p_d^i}^T, {p_a^j}^T)^T \in \mathbb{R}^4$, and the initial state of the system is defined as $x_0:=\left(({p_{d0}^i})^T,({p_{a0}^j})^T\right)^T=x(t_0)$. Here, $t_0$ represents the initial time of the game, and $t_f$ represents the termination time. The system's state equation is given by:
	\begin{equation}
		\dot x = Bu,
		\label{eq:sys}
	\end{equation}
	where
	$$
	B =\mathbf{I}_4 
	, \quad
	u = \begin{bmatrix} u_d^i \\ u_a^j \end{bmatrix},
	$$
	and $\mathbf{I}_4$ is the $4\times 4$ identity matrix.
	
The game is terminated if one of the following conditions is satisfied
	\begin{equation}
	\begin{aligned}
		&\text{1)}\quad \|p_d^i-p_a^j\|\leq R_d^c;
		\quad\text{2)}\quad \|p_a^j\|\leq R_{in};\\
		&\text{3)}\quad
		\begin{cases}
			\varphi_a^j<s_i\ \text{or}\ 
			\varphi_a^j\geq s_{i+1},
			& s_i<s_{i+1},\\
			s_{i+1}\leq\varphi_a^j<s_i,
			& s_i>s_{i+1}.
		\end{cases}
	\end{aligned}
	\label{eq:stop}
\end{equation}
	Here Condition 1 means the defender has captured the attacker; Condition 2 implies the attacker has entered the inner region of the target; Condition 3 indicates that the attacker has left the angular sector covered by defender $i$ and will be tackled by another defender. 
	
	The goals of the attacker and the defender are completely conflicting, making it suitable for analysis using zero-sum differential games. In such games, the strategies available to each agent depend on the information it has at the current time. Based on this, feedback control mappings can be constructed as follows:
	$$
	u_d^i(t) = \gamma_d^i(t,x(t)) , ~t\in [t_0,t_f],
	$$
	$$
	u_a^j(t) = \gamma_a^j(t,x(t)) , ~t\in [t_0,t_f],
	$$
	where $\gamma_d:[t_0,t_f]\times\mathbb{R}^4\longrightarrow\mathbb{R}^2$ and $\gamma_a:[t_0,t_f]\times\mathbb{R}^4\longrightarrow\mathbb{R}^2$. Consequently, the feedback strategy sets for the attacker and the defender are defined as follows:
	\begin{equation}
	\begin{aligned}
		\Gamma_d^i &:= \big\{\gamma_d^i \mid u_d^i(t)=\gamma_d(t,x(t)),\ t\in[t_0,t_f]\big\},\\
		\Gamma_a^j &:= \big\{\gamma_a^j \mid u_a^j(t)=\gamma_a(t,x(t)),\ t\in[t_0,t_f]\big\}.
	\end{aligned}
\end{equation}
	Define the following objective function:
	\begin{equation}
		\begin{aligned}
			J(u_d^i, u_a^j)
			=&\, {p_a^j}^T(t_f)F_a{p_a^j}(t_f) \\
			&-\big(p_d^i(t_f)-p_a^j(t_f)\big)^T
			F_d
			\big(p_d^i(t_f)-p_a^j(t_f)\big) \\
			&+\int_{t_0}^{t_f}
			\bigg[
			{p_a^j}^T(t)Q_a{p_a^j}(t) \\
			&\quad-\big(p_d^i(t)-p_a^j(t)\big)^T
			Q_d
			\big(p_d^i(t)-p_a^j(t)\big) \\
			&\quad+{u_a^j}^T(t)R_a{u_a^j}(t)
			-{u_d^i}^T(t)R_d{u_d^i}(t)
			\bigg]dt \\
			=&\,x^T(t_f)Fx(t_f) \\
			&+\int_{t_0}^{t_f}
			\left(
			x^T(t)Qx(t)+u^T(t)R_u u(t)
			\right)dt,
		\end{aligned}
		\label{eq:cost}
	\end{equation}
	where
	$F_a=f_a\otimes\mathbf{I}_2$,
	$Q_a=q_a\otimes\mathbf{I}_2$,
	$F_d=f_d\otimes\mathbf{I}_2$,
	$Q_d=q_d\otimes\mathbf{I}_2$,
	$R_a=r_a\otimes\mathbf{I}_2$, and
	$R_d=r_d\otimes\mathbf{I}_2$. Moreover,
	\[
	F=
	\begin{bmatrix}
		-F_d & F_d\\
		F_d & F_a-F_d
	\end{bmatrix},
	\quad
	Q=
	\begin{bmatrix}
		-Q_d & Q_d\\
		Q_d & Q_a-Q_d
	\end{bmatrix},
	\quad
	\]
	\[
	R_u=
	\begin{bmatrix}
		-R_d & \mathbf{0}\\
		\mathbf{0} & R_a
	\end{bmatrix}.
	\]
	Here, $f_a$, $f_d$, $q_a$, $q_d$, $r_a$, and $r_d$
	are positive tunable parameters, $\otimes$ denotes the
	Kronecker product, $\mathbf{I}_2$ represents the
	two-dimensional identity matrix, and $\mathbf{0}$ denotes
	the $2\times2$ zero matrix.
	
	The attacker aims to minimize the objective function, while the defender aims to maximize it. The feedback strategy mappings $({\gamma_d^i}^*,{\gamma_a^j}^*)$ constitute a feedback Nash equilibrium when the following inequality holds:
	\begin{equation}
		J\left( {\gamma_d^i}, {\gamma_a^j}^* \right)\le  J\left( {\gamma_d^i}^*, {\gamma_a^j}^* \right) \le J\left( {\gamma_d^i}^*, {\gamma_a^j} \right).
		\label{eq:nashi}
	\end{equation}
	In summary, the optimization formulation for the attacker and the defender are as follows:
	\begin{equation}
		\min_{\gamma_a^j} \max_{\gamma_d^i}~ J\left( {\gamma_d^i}, {\gamma_a^j} \right).
		\label{eq:game}
	\end{equation}
	In this game, the attacker and the defender continuously explore to find the best state feedback strategies to optimize~(\ref{eq:game}).

	\section{Theoretical Results}\label{sec:4.2}
	In this section, the role of the region partition algorithm in decomposing the attacker network is first introduced. Several key lemmas are presented to demonstrate that under the action of the region partition algorithm, each defender faces at most one attacker. Subsequently, the Nash equilibrium strategy set of the game is described using optimal control theory.
	
	\subsection{Key Lemmas}
	
	\begin{lemma}\label{le:se4}
		The region partition Equation~(\ref{eq:s}) ensures that the workload in each subregion is finally equal to the workload $2\frac{M\mathbf{e}(\pi)}{N}$.
	\end{lemma}
	\begin{IEEEproof}
		For a fixed attacker density $\rho(\theta)$, we construct the following Lyapunov function:
		$${V}(s) = \frac{1}{2}\sum_{i=1}^N\left(m_i-\overline{m}\right)^2,$$
		where $\overline{m}=\frac{1}{N}{\int_0^{2\pi}{\rho(\theta)d\theta}}$.
		The time derivative of the Lyapunov function along the trajectory of system~(\ref{eq:s}) is:
		\begin{equation*}
			\begin{aligned}
				\dot{V}(s) &= \sum_{i = 1}^N {\left( {{m_i} - \overline m} \right){{\dot m}_i}} \\
				&= \sum_{i = 1}^N {{m_i}{{\dot m}_i}} - \overline{m}\underbrace{\sum_{i=1}^{N}\dot{m}_i}_{0} \\
				&= \sum_{i=1}^N m_i \left( \rho(s_{i+1})\dot{s}_{i+1} - \rho(s_i)\dot{s}_i \right) \\
				&= \sum_{i=1}^N m_{i-1}\rho(s_i)\dot{s}_i - \sum_{i=1}^N m_i\rho(s_i)\dot{s}_i \\
				&= \sum_{i=1}^N \rho(s_i)\dot{s}_i (m_{i-1} - m_i) \\
				&= \sum_{i=1}^N \rho(s_i) \left[ \kappa_s(m_i - m_{i-1}) \right] (-(m_i - m_{i-1})) \\
				&= -\kappa_s\sum_{i=1}^N(m_i-m_{i-1})^2\rho\left({{s_i}}\right).
			\end{aligned}
		\end{equation*}
		Since $\kappa_s>0$ and $\rho \left({{s_i}}\right)>0$, it follows that $\dot{V}(s)\leq0$. Moreover, $\dot{V}(s)=0$ only when $m_i=m_{i-1}$ for all $i\in \mathbb{I}_N$. From Equation~(\ref{eq:hpi}), we have:
		$$
		\frac{1}{N}{\int_0^{2\pi}{\rho(\theta)d\theta}} = 2\frac{M\mathbf{e}(\pi)}{N} < \frac{M}{N}\sqrt{\frac{\pi}{h}}.
		$$
		Thus, the system's state approaches the equilibrium as $t\to\infty$, where $m_1(t)=m_2(t)=\cdots=m_N(t)= 2\frac{M\mathbf{e}(\pi)}{N}$.
	\end{IEEEproof}
	
	\begin{lemma}
		There exist positive constants $c_1$ and $c_2$ such that for system~(\ref{eq:s}), $|m_i(t)-m_{i-1}(t)|\leq c_1 e^{-c_2t}$ for all $t\geq0$.
	\end{lemma}
	\begin{IEEEproof}
		The proof can be found in Lemma 3.2 and Lemma 3.3 of~\cite{ZHAI2023111246}.
	\end{IEEEproof}
	
	Before proving that the region partition ensures each defender faces at most one attacker, we first need to obtain the minimum absolute phase difference between attackers. The minimum absolute phase difference $\Delta\varphi_{\min}$ between attackers can be obtained by solving:
	\begin{equation}
			 \Delta\varphi_{\min} = \min_{ i,j\in I_M, i\neq j}|\delta(\varphi_a^i,\varphi_a^j)|	
		\label{eq:minpa}
	\end{equation}
	
	Based on Equation~(\ref{eq:minpa}), the minimum absolute phase difference $\Delta\varphi_{\min}$ between attackers can be used to establish a relationship between $h$, $\Delta\varphi_{\min}$, $M$, and $N$, ensuring that the region partition algorithm guarantees at most one attacker per defender's subregion. Let $N^* > 0$ be a scaling constant uniquely determined by the system load ratio $M/N$. Specifically, $N^*$ represents the standard score such that the cumulative probability of a normal distribution $\mathcal{N}\left(0,\frac{1}{2h}\right)$ within the interval $[-N^*, N^*]$ equals the ratio $M/N$. This constant can be readily obtained from standard normal distribution tables.
	
	\begin{theorem}
		For any $\Delta\varphi_{\min}>0$, if the scaling factor $h$ satisfies $h\ge\frac{{N^*}^2}{2\Delta\varphi_{\min}^2}$, then, at every equilibrium of the region-partition dynamics~(\ref{eq:s}),
each defender's subregion contains at most one attacker.
		\label{se4:th1}
	\end{theorem}
	\begin{IEEEproof}
		By way of contradiction, assume that there exists a defender's subregion with more than one attacker. Let the defender be denoted by index $i$, and there are at least two attackers $j$ and $j+1$ within the subregion. Then, according to Equation~(\ref{eq:bi_m}), we have
		\begin{equation*}
			\begin{aligned}
				m_i &= \int_{E_i} \rho(\theta) d\theta \\
				&\ge \int_{\varphi_a^j}^{\varphi_a^{j+1}} \sum_{j\in I_M^i} K(\theta,\varphi_a^j,h)d\theta \\
				&\ge 2 \int_0^{\Delta\varphi_{\min}} e^{-h\eta^2}d\eta = 2\mathbf{e}(\Delta\varphi_{\min}).
			\end{aligned}
		\end{equation*}
		The first inequality holds because of equations~(\ref{eq:rho8}) and~(\ref{eq:rhoi}), the second inequality holds because the phase range $E_i$ must include attackers $j$ and $j+1$, and the third inequality is valid because the phase difference between attackers must be greater than the minimum absolute phase difference.
		
		By Lemma~\ref{le:se4}, $m_i=2\frac{M\mathbf{e}(\pi)}{N}$. Hence, the following inequality must hold:
		\begin{equation*}
			\begin{aligned}
				2\frac{M\mathbf{e}(\pi)}{N} > 2\mathbf{e}(\Delta\varphi_{\min}),
			\end{aligned}
		\end{equation*}
		which simplifies to
		$$
		\frac{M}{N} > \frac{\mathbf{e}(\Delta\varphi_{\min})}{\mathbf{e}(\pi)}.
		$$
		Observing the right-hand side of the inequality, we have
		\begin{equation*}
			\begin{aligned}
				\frac{\mathbf{e}(\Delta\varphi_{\min})}{\mathbf{e}(\pi)} &> \frac{\mathbf{e}(\Delta\varphi_{\min})}{\mathbf{e}(\infty)}\\
				&= \frac{\int_0^{\Delta\varphi_{\min}} e^{-h\eta^2}d\eta}{\int_0^{+\infty} e^{-h\eta^2}d\eta}.
			\end{aligned}
		\end{equation*}
		Observing the last term, we notice that it can be transformed into the probability of the values falling within the interval $(-{\Delta\varphi_{\min}},+{\Delta\varphi_{\min}})$ for a normal distribution $X\sim N(0,\sqrt{\frac{1}{2h}})$ with mean 0 and standard deviation $\sqrt{\frac{1}{2h}}$, denoted as $\Phi(\Delta\varphi_{\min})$. Thus, the inequality simplifies to
		\begin{equation*}
			\frac{M}{N}> \Phi(\Delta\varphi_{\min}).
		\end{equation*}
		Since $ \frac{M}{N} \in (0,1)$, $\frac{M}{N}$ can also be represented using the probability $\Phi(\cdot)$. Let $\frac{M}{N}=\Phi(N^*\sqrt{\frac{1}{2h}})$, where $N^*$ is a positive number. Therefore,
		\begin{equation*}
			\Phi(N^*\sqrt{\frac{1}{2h}})> \Phi(\Delta\varphi_{\min}).
		\end{equation*}
		Since $\Phi(\cdot)$ is monotonically increasing, we have $N^*\sqrt{\frac{1}{2h}}>\Delta\varphi_{\min}$, which simplifies to $h<\frac{{N^*}^2}{2\Delta\varphi_{\min}^2}$. Since $h\ge\frac{{N^*}^2}{2\Delta\varphi_{\min}^2}$, it contradicts the original assumption. Therefore, the original proposition does not hold, implying that there is at most one attacker within each defender's subregion in the region partition dynamics~(\ref{eq:s}).
	\end{IEEEproof}

	\begin{remark}
		In statistics, the well-known ``68-95-99.7 rule" exists for normal distributions, which is commonly used to describe the characteristics of data distributions~\cite{moore2017introduction}. In the context of normal distribution, approximately 68\% of data points fall within one standard deviation of the mean, while 95\% of data points are distributed within two standard deviations of the mean. Approximately 99.7\% of data points are distributed within three standard deviations of the mean. Hence, by looking up the standard normal distribution table, we can determine the number of standard deviations corresponding to data within a specific percentage. In the standard normal distribution table, the standard normal distribution scores corresponding to different percentile points are usually provided. Therefore, $N^*$ can be obtained by looking up the table.
	\end{remark}
	\begin{remark}
	When multiple attackers have identical phase angles, i.e.,
	$\Delta\varphi_{\min}=0$, Theorem~\ref{se4:th1} is no longer
	applicable. Under the considered radial-approach geometry, if
	these attackers remain phase-aligned while moving toward the
	target, they follow the same interception path. Therefore, the
	defender can first intercept the attacker closest to the target
	and subsequently engage the following attackers along the same
	path. The allocation property for this overlapping-phase case
	is further analyzed in Lemma~\ref{le:NKk} and
	Theorem~\ref{se4:th2}.
\end{remark}

	Given that the number of defenders exceeds that of attackers ($N > M$), certain regions will inevitably contain no attackers. Defenders assigned to these regions are referred to as ``idle'' agents, denoted by the set $I_d^r$. Although not currently engaged in interception, these agents must be strategically positioned to guard against potential future threats.

    For agents executing the coverage strategy, the optimal standby position $p_i^*$ is determined by calculating the density-weighted centroid of the sector $E_i$. Let $q$ be an arbitrary point within the sector $E_i$, and let $\theta_q$ denote the phase angle of $q$. We define the spatial density function $\tilde{\rho}_i(q)$ by mapping the angular density to the planar domain based on the presence of attackers:

\begin{equation}
    \tilde{\rho}_i(q) = 
    \begin{cases}
        \rho_i^*(\theta_q), & \text{if } I_M^i \neq \emptyset \\
        \rho_i(\theta_q), & \text{if } I_M^i = \emptyset
    \end{cases}
    \label{eq:switching_rho}
\end{equation}
where $\rho_i^*(\cdot)$ and $\rho_i(\cdot)$ are the local and global density functions defined in Equation~(\ref{eq:rho8}) and Equation~(\ref{eq:rhoi}), respectively. The optimal position $p_i^*$ is then obtained by minimizing the weighted squared Euclidean distance, i.e., the expected value of a density-weighted cost function over the sector $E_i$:
\begin{equation}
    p_i^* = \arg \min_{p \in E_i} \int_{E_i} f(p, q) \tilde{\rho}_i(q) dq.
    \label{eq:se4H}
\end{equation}
where $f(p, q)$ represents a distance-dependent cost function measuring the cost between the defender's position $p$ and a point $q$ in the domain. This formulation ensures that when attackers are detected within the subregion ($I_M^i \neq \emptyset$), the centroid is calculated solely based on the local density $\rho_i^*$ to maintain precise tracking. Conversely, when the subregion is empty ($I_M^i = \emptyset$), the defender utilizes the aggregated global density $\rho_i$ to move towards high-risk boundaries for early warning.

	To drive the idle defender $i$ toward this optimal position, the following control input is designed:
	\begin{equation}
		u_d^i = -\kappa_d (p_d^i - p_i^*).
		\label{eq:idr_u}
	\end{equation}
	Here, $\kappa_d$ is a positive control gain, ensuring that the defender converges to the global optimal solution within the 2D plane.
	
	It is observed that unassigned (``idle'') agents can be strategically positioned to cover potential risks, whereas agents already engaged in the game lack the capacity to handle additional threats. The following theoretical analysis validates the effectiveness of this allocation mechanism. According to Theorem~\ref{se4:th1}, provided that the phases of attackers do not overlap, the region partitioning algorithm ensures that each defender confronts at most one attacker. Consequently, even if an attacker traverses into a neighboring region, the partition dynamics will adapt to redistribute the workload, thereby maintaining the one-to-one engagement constraint and preventing any defender from facing multiple attackers simultaneously.
	
	\begin{figure}
		\centering
		\includegraphics[width=0.45\textwidth]{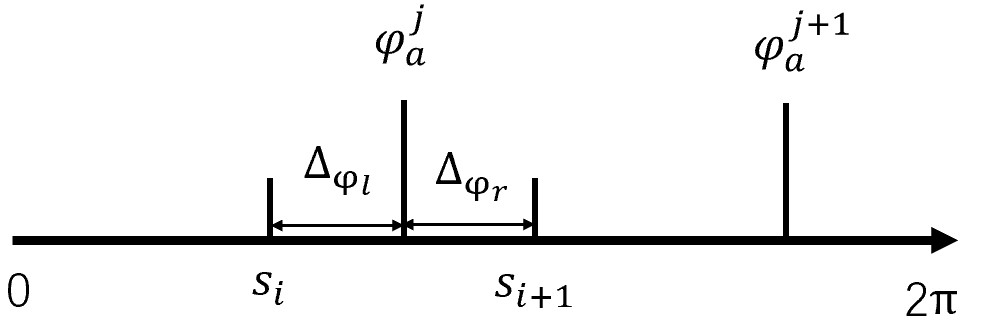}
		\caption{Phase relationship used in Lemma~\ref{le:NKk}: overlapping attacker phases and their angular distances to adjacent partition boundaries.}
		\label{fig:lema}
	\end{figure}
	
	\begin{lemma}
		When $\Delta\varphi_{\min}=0$, i.e., $k$ attackers have overlapping phases ($k$ is a positive integer, $k>1$), the minimum distance from the attackers to the boundary satisfies the following inequality:
		$$
		0\le \Delta_\varphi < N_k^* \sqrt{\frac{1}{2h}}.
		$$
		Here, $\Delta_\varphi = \min \{\Delta_{\varphi_r},\Delta_{\varphi_l}\}$, $\Delta_{\varphi_r}$ denotes the absolute value of the phase difference between an attacker and the partition angle $s_{i+1}$, and $\Delta_{\varphi_l}$ denotes the phase difference between an attacker and the partition angle $s_i$ (see Figure~\ref{fig:lema}). $N_k^*$ satisfies $\Phi(N_k^* \sqrt{\frac{1}{2h}} ) = \frac{M}{kN}$.
		\label{le:NKk}
	\end{lemma}
	\begin{IEEEproof}
		By Lemma~\ref{le:se4}, it is known that the workload of defender $i$ satisfies the following inequality:
		$$
		k\left({\mathbf{e}(\Delta_{\varphi_r})}+{\mathbf{e}(\Delta_{\varphi_l})}\right) \le \frac{2M}{N} {\mathbf{e}(\pi)}.
		$$
		Equality holds only when $k=M$.
		Since $\Delta_\varphi = \min \{\Delta_{\varphi_r},\Delta_{\varphi_l}\}$, we have
		$$
		\left({\mathbf{e}(\Delta_{\varphi_r})}+{\mathbf{e}(\Delta_{\varphi_l})}\right)\ge2{\mathbf{e}(\Delta_{\varphi})}.
		$$
		Equality holds only when $\Delta_{\varphi_r}=\Delta_{\varphi_l}$. Substituting into the above inequality, we obtain
		$$
		{\mathbf{e}(\Delta_{\varphi})}\le \frac{M}{kN} {\mathbf{e}(\pi)}.
		$$
		Since $\mathbf{e}(\pi)<\mathbf{e}(\infty)$, we have
		\[
		\mathbf{e}(\Delta_\varphi)
		\leq \frac{M}{kN}\mathbf{e}(\pi)
		<\frac{M}{kN}\mathbf{e}(\infty).
		\]
		Therefore,
		\[
		\Phi(\Delta_\varphi)
		<
		\Phi\left(N_k^*\sqrt{\frac{1}{2h}}\right).
		\]
		Since $\Phi(\cdot)$ is monotonically increasing, we have $\Delta_{\varphi} < N_k^*\sqrt{\frac{1}{2h}}$.
	\end{IEEEproof}
	
	Based on Lemma~\ref{le:NKk}, the following theorem can be derived:
	
	\begin{theorem}
	\label{se4:th2}
	When $k$ attackers have overlapping phases, where $2\le k<M$, if the scaling factor $h$ satisfies $h\geq\frac{{N^*}^2}{2(\Delta\varphi_{\min}^*-\Delta_\varphi)^2}$, then, at every equilibrium of the region-partition dynamics~(\ref{eq:s}), at least one of the neighboring defenders of defender agent $i$ is an ``idle'' defender. Here, $\Delta\varphi_{\min}^*$ denotes the angular distance between the overlapping attackers and the nearest non-overlapping attacker on the same side as the closest adjacent partition boundary, and $N^*$ is the same as in Theorem~\ref{se4:th1}.
\end{theorem}
	
	\begin{IEEEproof}
	The above theorem can be reformulated as follows:
	There exists $\psi\in\{j-1,j+1\}$ such that
	\[
	\int_{s_{i+1}}^{\varphi_a^\psi}\rho(\theta)d\theta
	>
	2\frac{M\mathbf{e}(\pi)}{N}.
	\]
	Due to Lemma~\ref{le:se4}, the above inequality implies that
	there must exist an ``idle'' defender within the angular interval
	defined by $[s_{i+1},\varphi_a^\psi]$.

	Since $\Delta\varphi_{\min}^*$ denotes the angular distance
	between the overlapping attackers and the nearest non-overlapping
	attacker on the same side as the closest adjacent partition
	boundary, the following inequality holds:
	\[
	\int_{s_{i+1}}^{\varphi_a^\psi}\rho(\theta)d\theta
	>
	\mathbf{e}(\Delta\varphi_{\min}^*-\Delta_\varphi)
	+k(\mathbf{e}(\Delta\varphi_{\min}^*)-\mathbf{e}(\Delta_\varphi)).
	\]
	Therefore, it is sufficient to prove that
	\begin{equation}
	\begin{aligned}
	&\frac{\mathbf{e}(\Delta\varphi_{\min}^*-\Delta_\varphi)}
	{\mathbf{e}(\pi)}
	+k\left(
	\frac{\mathbf{e}(\Delta\varphi_{\min}^*)}
	{\mathbf{e}(\pi)}
	-\frac{\mathbf{e}(\Delta_\varphi)}
	{\mathbf{e}(\pi)}
	\right)
	>
	2\frac{M}{N}.
	\end{aligned}
	\label{eq:the1}
	\end{equation}

	According to the definition of $N^*$ and the condition of $h$,
	we have
	\[
	\frac{\mathbf{e}(\Delta\varphi_{\min}^*-\Delta_\varphi)}
	{\mathbf{e}(\infty)}
	\geq
	\frac{M}{N}.
	\]
	Since $\mathbf{e}(\pi)<\mathbf{e}(\infty)$, it follows that
	\[
	\frac{\mathbf{e}(\Delta\varphi_{\min}^*-\Delta_\varphi)}
	{\mathbf{e}(\pi)}
	>
	\frac{M}{N}.
	\]
	Moreover, because $\mathbf{e}(\cdot)$ is monotonically increasing,
	we have
	\[
	\frac{\mathbf{e}(\Delta\varphi_{\min}^*)}
	{\mathbf{e}(\pi)}
	>
	\frac{M}{N}.
	\]
	According to Lemma~\ref{le:NKk}, we obtain
	\[
	\frac{\mathbf{e}(\Delta_\varphi)}
	{\mathbf{e}(\pi)}
	\leq
	\frac{M}{kN}.
	\]
	Therefore, the left-hand side of inequality~(\ref{eq:the1})
	satisfies
	\[
	\begin{aligned}
	&\frac{\mathbf{e}(\Delta\varphi_{\min}^*-\Delta_\varphi)}
	{\mathbf{e}(\pi)}
	+k\left(
	\frac{\mathbf{e}(\Delta\varphi_{\min}^*)}
	{\mathbf{e}(\pi)}
	-\frac{\mathbf{e}(\Delta_\varphi)}
	{\mathbf{e}(\pi)}
	\right)\\
	&>
	\frac{M}{N}
	+k\left(
	\frac{M}{N}
	-\frac{M}{kN}
	\right)
	=
	k\frac{M}{N}
	\geq
	2\frac{M}{N}.
	\end{aligned}
	\]
	Hence, inequality~(\ref{eq:the1}) holds, which yields
\[
\int_{s_{i+1}}^{\varphi_a^\psi}\rho(\theta)d\theta
>2\frac{M\mathbf{e}(\pi)}{N}.
\]
It follows that $s_{i+2}<\varphi_a^\psi$. Since no attacker
lies between $s_{i+1}$ and $\varphi_a^\psi$, the subregion
$E_{i+1}$ contains no attacker. Therefore, defender agent
$i+1$ is an ``idle'' defender. This completes the proof.
\end{IEEEproof}
	
	\begin{remark} Theorem~\ref{se4:th2} indicates that when there are attackers having  the same phase angles, a defender may be responsible for two attackers. But there must be an ``idle'' defender among its neighbors. When attackers no longer have overlapping phases, according to Theorem~\ref{se4:th1}, all agents will be assigned to at most one attacker, which means attackers will be assigned to the ``idle'' defender.
	\end{remark}
	
	\subsection{Game Analysis}
	Building upon the decoupling mechanism established in the previous section, we now transition to the analysis of the localized one-to-one barrier interception game involving a single attacker and a single defender. For notational simplicity, the attacker and defender indices are omitted hereafter, and Equation~(\ref{eq:sys}) can be decomposed as follows:
	\begin{equation}
		\dot x = B_du_d+B_au_a,
		\label{eq:xbdba}
	\end{equation}
	where
	$$
	B_d =
	\begin{bmatrix}
		1 & 0\\
		0 & 1 \\
		0 & 0 \\
		0 & 0
	\end{bmatrix}
	,\quad
	B_a =
	\begin{bmatrix}
		0 & 0\\
		0 & 0 \\
		1 & 0\\
		0 & 1
	\end{bmatrix}.
	$$
	There have been numerous research results on solving Nash equilibrium strategies for two-player zero-sum differential games. By leveraging the theoretical framework established in~\cite{engwerda2005lq} and integrating it with our model, we provide rigorous proofs for the existence of Nash equilibrium and Nash equilibrium strategies.
	
	\begin{lemma}
		Consider a zero-sum differential game~(\ref{eq:game}) with one attacker and one defender. If the following Riccati differential equation
		\begin{equation}
			\dot{P}=PB_aR_a^{-1}B_a^TP-PB_dR_d^{-1}B_d^TP-Q ,P(t_f) = F
			\label{eq:dotp}
		\end{equation}
		has a symmetric solution $P(t)$ and the strategy set $(\gamma_d,\gamma_a)$ satisfies $\gamma_d = R_d^{-1}B_d^TP(t)x(t)$ and $\gamma_a = -R_a^{-1}B_a^TP(t)x(t)$, then the zero-sum differential game~(\ref{eq:game}) reaches a feedback Nash equilibrium.
	\end{lemma}
	
	\begin{IEEEproof}
		Define the matrix $P$, then we can handle the cost function~(\ref{eq:cost}) as follows:
		\begin{equation*}
			\begin{aligned}
				J(u_d, u_a) &= x^T(t_f) F x(t_f) + \int_{t_0}^{t_f} \left( x^TQx + u^TR_u u\right) dt \\
				&\quad + \int_{t_0}^{t_f} \frac{d}{dt}(x^TPx) dt - \int_{t_0}^{t_f} \frac{d}{dt}(x^TPx) dt \\
				&= x^T(t_0) P(t_0) x(t_0) \\
				&\quad + \int_{t_0}^{t_f} \left( x^TQx + u^TR_u u +\frac{d}{dt}(x^TPx) \right) dt
			\end{aligned}
		\end{equation*}
		where $P(t_f) = F$.
		Consider the integral term, we get
		{\small
		\begin{equation*}
			\begin{aligned}
				& x^TQx + u^TR_u u +\frac{d}{dt}(x^TPx) \\
				&\quad = x^TQx + u^TR_u u + \dot{x}^TPx + x^T\dot{P}x + x^TP\dot{x} \\
				&\quad = x^TQx + u_a^TR_a u_a - u_d^T R_d u_d \\
				&\qquad + (B_du_d+B_au_a)^TPx + x^T\dot{P}x + x^TP(B_du_d+B_au_a) \\
				&\quad = x^T(Q + \dot{P})x + u_a^TR_au_a - u_d^TR_du_d \\
				&\qquad + u_d^TB_d^TPx + u_a^TB_a^TPx + x^TPB_du_d + x^TPB_au_a \\
				&\quad = (u_a+R_a^{-1}B_a^TPx)^T R_a (u_a+R_a^{-1}B_a^TPx) \\
				&\qquad - (u_d-R_d^{-1}B_d^TPx)^T R_d (u_d-R_d^{-1}B_d^TPx) \\
				&\qquad + x^T\big(Q + \dot{P} - PB_aR_a^{-1}B_a^TP \\
				&\qquad\qquad + PB_dR_d^{-1}B_d^TP\big)x
			\end{aligned}
		\end{equation*}
		}
		Setting $\dot{P}=PB_aR_a^{-1}B_a^TP-PB_dR_d^{-1}B_d^TP-Q$ yields
		\begin{equation*}
			\begin{aligned}
				J &= x^T(t_0) P(t_0) x(t_0) \\
				&\quad + \int_{t_0}^{t_f} \bigg[ (u_a + R_a^{-1}B_a^T P x)^T R_a (u_a + R_a^{-1}B_a^T P x) \\
				&\qquad - (u_d - R_d^{-1}B_d^T P x)^T R_d (u_d - R_d^{-1}B_d^T P x) \bigg] dt
			\end{aligned}
		\end{equation*}
		Let $u_d^*=R_d^{-1}B_d^TPx$ and
		$u_a^*=-R_a^{-1}B_a^TPx$. Then,
		\begin{equation*}
			\begin{aligned}
				J(u_d,u_a^*)
				={}&x^T(t_0)P(t_0)x(t_0)\\
				&-\int_{t_0}^{t_f}
				(u_d-u_d^*)^TR_d(u_d-u_d^*)dt\\
				\leq{}&J(u_d^*,u_a^*),
			\end{aligned}
		\end{equation*}
		and
		\begin{equation*}
			\begin{aligned}
				J(u_d^*,u_a)
				={}&x^T(t_0)P(t_0)x(t_0)\\
				&+\int_{t_0}^{t_f}
				(u_a-u_a^*)^TR_a(u_a-u_a^*)dt\\
				\geq{}&J(u_d^*,u_a^*).
			\end{aligned}
		\end{equation*}
		Therefore, the saddle-point condition~(\ref{eq:nashi}) holds.

	\end{IEEEproof}
	
	Since the matrices $B_d$ and $B_a$ are constant, the Riccati Equation~(\ref{eq:dotp}) is a standard matrix differential equation. By solving the above Riccati differential equation, Nash equilibrium strategies can be obtained:
	\begin{equation}
		u_d(t) = R_d^{-1}B_d^TP(t)x(t),\quad u_a(t)= -R_a^{-1}B_a^TP(t)x(t).
		\label{eq:uaud}
	\end{equation}
	Equation~(\ref{eq:uaud}) provides the nominal Nash equilibrium
strategies for the unconstrained local game. In the numerical
implementation, these strategies are further modified to satisfy
the velocity and boundary constraints.

	\begin{remark}
		The control strategy described in Equation~(\ref{eq:uaud}) is consistent with the control strategy derived from $H_\infty$ control, and the relationship between these formulations has been established in the literature~\cite{basar2008hinfinity}.
	\end{remark}

\subsection{Algorithm Design}

The distributed barrier interception algorithm operates with an update period $T$ and consists of three sequential phases: information fusion, region partitioning, and defense control.

In the information fusion phase, the defense network acquires the phase distribution of attackers within the detection zone $R_{in} \le \|p_a\| \le R_{detect}$. The scaling factor $h$ of the Gaussian kernel is adaptively adjusted based on the minimum phase difference $\Delta\varphi_{\min}$ among attackers. To mitigate the slow convergence caused by an excessively large
scaling factor, a saturation constraint is applied to $h$, defined as
$h=\min\left\{\frac{{N^*}^2}{2\Delta\varphi_{\min}^2},h_{\max}\right\}$,
where $h_{\max}$ is a preset threshold. When the theoretical value
does not exceed $h_{\max}$, the condition in
Theorem~\ref{se4:th1} is satisfied. Otherwise, the attackers are
nearly overlapping in phase, and $h$ is bounded to avoid the slow
convergence caused by an excessively narrow Gaussian kernel. In this
case, the strict allocation guarantee in Theorem~\ref{se4:th1} is
relaxed in favor of numerical convergence. Consequently, the global density function $\rho(\theta)$ is synthesized using the modified $h$.

In the region partitioning phase, each defender node iteratively updates the partition boundaries $s_i$ based on load balancing principles. A topology preservation mechanism is incorporated to maintain assignment stability: the update of a boundary $s_i$ is executed only if its angular distance to any attacker within the engagement zone $\|p_a\| \le R_{out}$ exceeds a threshold $\epsilon$. This process iterates $K^*$ times per cycle to ensure boundary convergence, thereby defining the responsibility region $E_i$ for each defender.

In the control execution phase, defender $i$ switches its control mode based on the state of its assigned region $E_i$. If an attacker $j$ exists within the region and satisfies
$\|p_a^j\|\leq R_{out}$, the defender solves the Riccati equation
and computes the nominal game-based control input according to
Equation~(\ref{eq:uaud}). The control input is then modified by
velocity saturation and a boundary safety correction before
execution. Conversely, if the region is free of immediate threats,
the defender executes the coverage strategy. In this mode, the algorithm selects either the local density $\rho_i^*$ or the global density $\rho_i$ as the weighting function based on the presence of distant targets, driving the defender toward the weighted centroid.

\begin{algorithm}[t]
    \caption{\label{tab:allan}Barrier Interception Algorithm with Adaptive Target Assignment}
    \hspace*{0.02in} \textbf{Initialize:} $m=0, T, K^*, R_{out}, \epsilon, h_{\max}$\\
    \hspace*{0.02in} For any defender $i \in I_d$, execute at each time step:
    \begin{algorithmic}[1]
        \WHILE{True}
            \IF{$mod(m,T) = 0$}
                \STATE Obtain attacker phases $\varphi_a$ via Algorithm~\ref{tab:algong}
                \STATE Compute $\Delta\varphi_{\min} = \min_{i \neq j} |\delta(\varphi_a^i, \varphi_a^j)|$
                \STATE Update kernel width: $h \gets \min\left(\frac{{N^*}^2}{2\Delta\varphi_{\min}^2}, h_{\max}\right)$
                \STATE Construct global density $\rho(\theta)$ using updated $h$
                
                \FOR{$k \gets 1$ to $K^*$}
                    \STATE Compute load $m_i$ and boundary velocity $\dot{s}_i$ via~(\ref{eq:s})
                    \IF{$|\delta(s_i, \varphi_a^j)| > \epsilon, \forall j \text{ s.t. } R_{in} \le \|p_a^j\| \le R_{out}$}
                        \STATE $s_i \gets s_i + \dot{s}_i \cdot \Delta t$
                    \ENDIF
                \ENDFOR
            \ENDIF
            
            \STATE Identify local attacker set $I_M^i$ within $E_i = [s_i, s_{i+1})$
            
            \IF{$\exists j \in I_M^i \text{ s.t. } \|p_a^j\| \le R_{out}$}
                \STATE Solve the Riccati equation for $P$
				\STATE Compute nominal game-based controls via Equation~(\ref{eq:uaud})
				\STATE Apply velocity saturation and boundary safety correction
            \ELSE
                \IF{$I_M^i \neq \emptyset$}
                    \STATE $\tilde{\rho} \gets \rho_i^*$ 
                \ELSE
                    \STATE $\tilde{\rho} \gets \rho_i$ 
                \ENDIF
                \STATE Compute centroid $p_i^*$ using $\tilde{\rho}$ via Equation~(\ref{eq:centroid_calc})
                \STATE Update $u_d^i$ via Coverage Control Equation~(\ref{eq:idr_u})
            \ENDIF
            \STATE $m \gets m+1$
        \ENDWHILE
    \end{algorithmic}
\end{algorithm}
	
\section{Numerical Simulation} \label{sec:4.3}

This section evaluates the effectiveness, scalability, and computation time of the proposed distributed interception algorithm through numerical simulations. The validation is structured into two primary components: (1) a detailed analysis of the effectiveness of the pre-positioning mechanism and the scalability of the algorithm, and (2) a comparative analysis against a centralized strategy to highlight differences in operational characteristics and computational efficiency.

The control laws in Equation~(\ref{eq:uaud}) are used as nominal
game-based commands in the simulations. To account for the physical
speed limits, the attacker and defender commands are saturated at
$U_{\max}^a$ and $U_{\max}^d$, respectively. Moreover, a boundary
safety correction is applied to the defender command to prevent
the defender from entering the protected region
$\|p_d^i\|<R_{in}$. A projection onto the inner boundary is used
only to eliminate residual constraint violations caused by the
discrete-time numerical integration. Therefore, the implemented
commands constitute a constrained numerical realization of the
nominal strategies rather than exact Nash equilibrium strategies.

Throughout the simulations, the coverage cost function in Equation~(\ref{eq:se4H}) is chosen as the squared Euclidean distance, i.e., $f(p, q) = \|p - q\|^2$. Accordingly, the optimal standby position $p_i^*$ is given by the density-weighted centroid of sub-region $E_i$:
\begin{equation}
    p_i^* = \frac{\int_{E_i} q \tilde{\rho}_i(q) dq}{\int_{E_i} \tilde{\rho}_i(q) dq}.
    \label{eq:centroid_calc}
\end{equation}

\begin{table}[h]
\centering
\caption{Simulation Parameters}
\label{tab:sim_params}
\footnotesize
\setlength{\tabcolsep}{3pt}
\begin{tabular}{l c p{3.8cm}}
\hline
Parameter & Value & Description \\
\hline
$N, M$ & $9, 8$ & Numbers of defenders and attackers \\
$R_{in}, R_{out}$ & $20, 50~\text{m}$ & Protection and engagement radii \\
$R_{detect}, R_d^c$ & $100, 0.5~\text{m}$ & Detection and capture radii \\
$U_{\max}^d, U_{\max}^a$ & $4.0, 5.0~\text{m/s}$ & Maximum speeds \\
$\kappa_s, \kappa_d$ & $0.02, 5.0$ & Control gains \\
$h_{max}$ & $60.0$ & Kernel width limit \\
$q_d, q_a, r_d, r_a$ & $40.0, 1.0, 0.01, 0.5$ & Game weights \\
\hline
\end{tabular}
\end{table}

\subsection{Validation of Effectiveness and Scalability}

The first part of the simulation validates the effectiveness of the proposed algorithm and its scalability across varying swarm sizes. A scenario with $N=9$ defenders and $M=8$ attackers is considered to examine the two-stage defense mechanism. The defense area is defined by a protection boundary $R_{in} = 20~\text{m}$ and an engagement boundary $R_{out} = 50~\text{m}$, with attackers initialized at $R_{detect} = 100~\text{m}$. To evaluate the algorithm under unfavorable kinematic conditions, the maximum speed of defenders is restricted to $U_{\max}^d = 4.0~\text{m/s}$, while attackers move at $U_{\max}^a = 5.0~\text{m/s}$. Capture is defined to occur when the relative distance is within $R_{d}^c = 0.5~\text{m}$. Detailed simulation parameters are summarized in Table~\ref{tab:sim_params}.

\begin{figure}[h]
    \centering
    \includegraphics[width=0.92\linewidth]{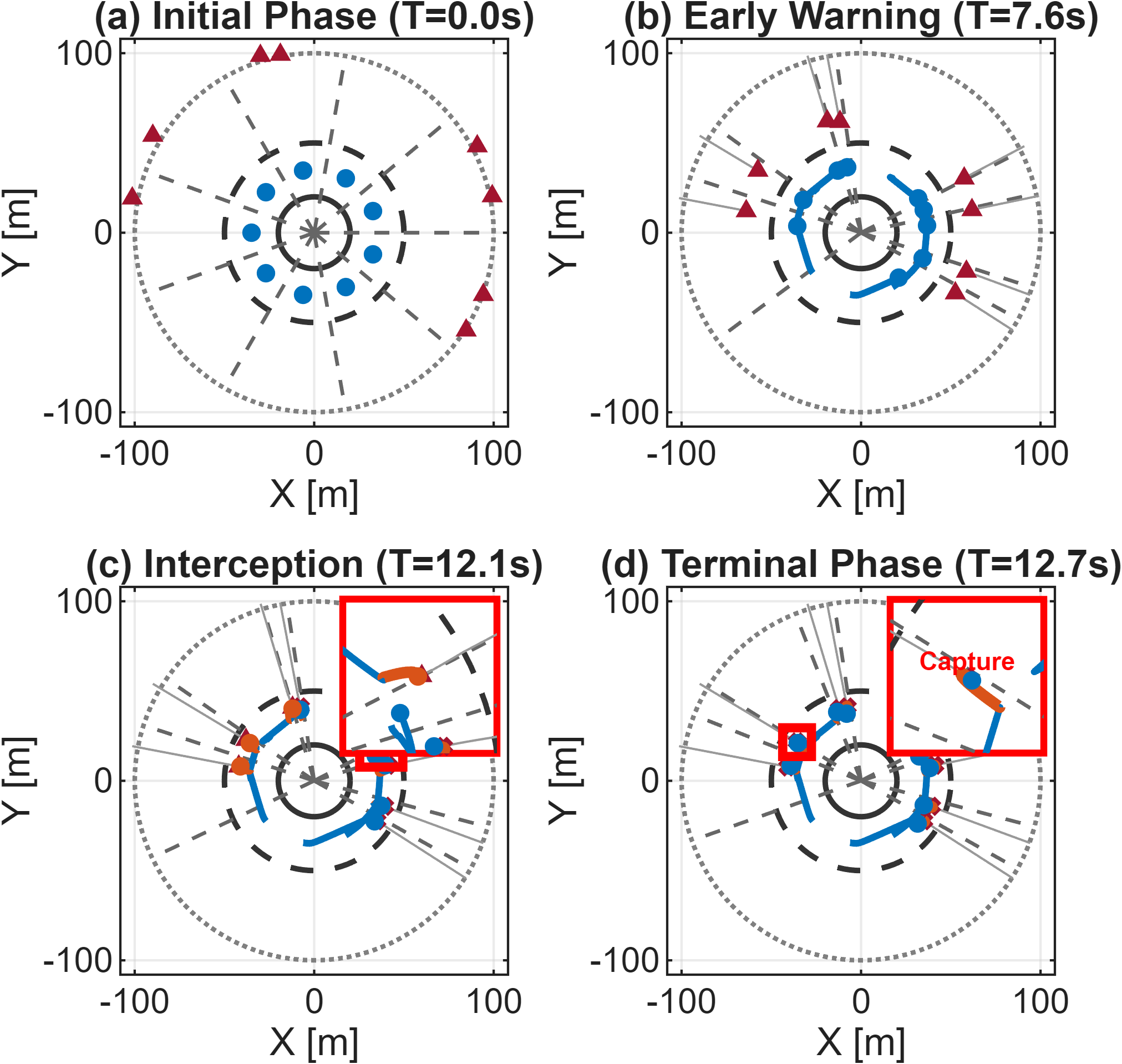}
    \caption{Snapshots of Feasibility Verification. (a) Initial configuration; (b) Defenders pre-position to centroids; (c) Mode switching (Blue $\to$ Red); (d) Capture of faster targets.}
    \label{fig:exp1_snapshots}
\end{figure}

The simulation snapshots in Fig.~\ref{fig:exp1_snapshots} illustrate the effectiveness of the pre-positioning mechanism. During the early-warning phase ($t = 7.60~\text{s}$), defenders maneuver toward the density-weighted centroids of their assigned regions. This proactive deployment establishes a geometric advantage at the engagement boundary before attackers penetrate the defense zone. Once the attackers enter the engagement zone, defenders transition to the interception mode, and all attackers are successfully intercepted in this scenario. The results demonstrate that the algorithm establishes an advantageous interception geometry, thereby enabling successful interception even when the defenders are slower than the attackers.

\begin{figure}[h]
    \centering
    \begin{subfigure}[b]{0.47\columnwidth}
        \centering
        \includegraphics[width=\textwidth]{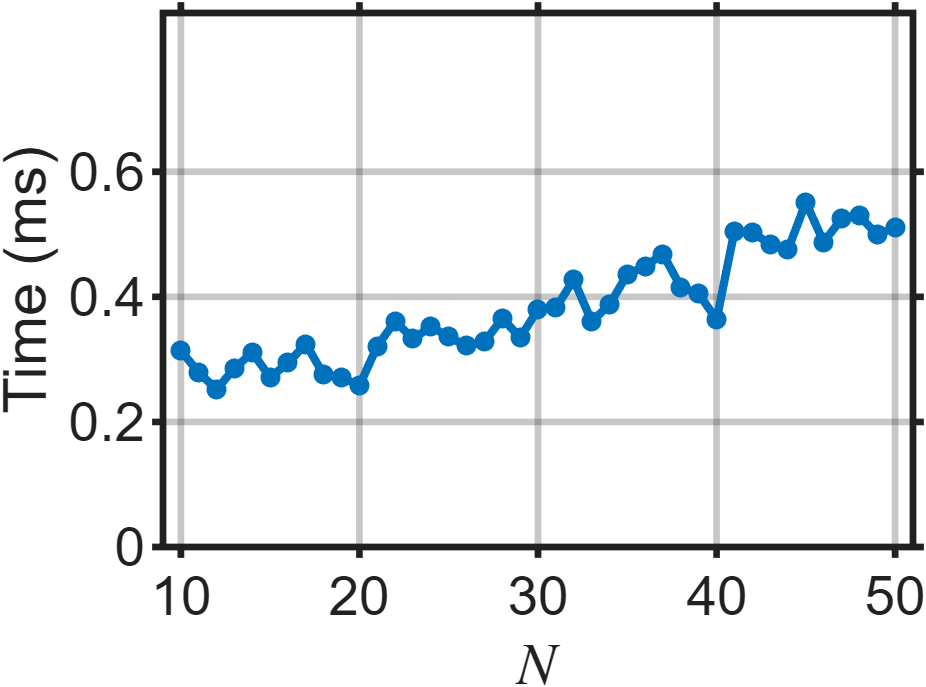}
        \caption{Computation time versus system scale $N$}
        \label{fig:sim2_scalability}
    \end{subfigure}
    \hfill
    \begin{subfigure}[b]{0.47\columnwidth}
        \centering
        \includegraphics[width=\textwidth]{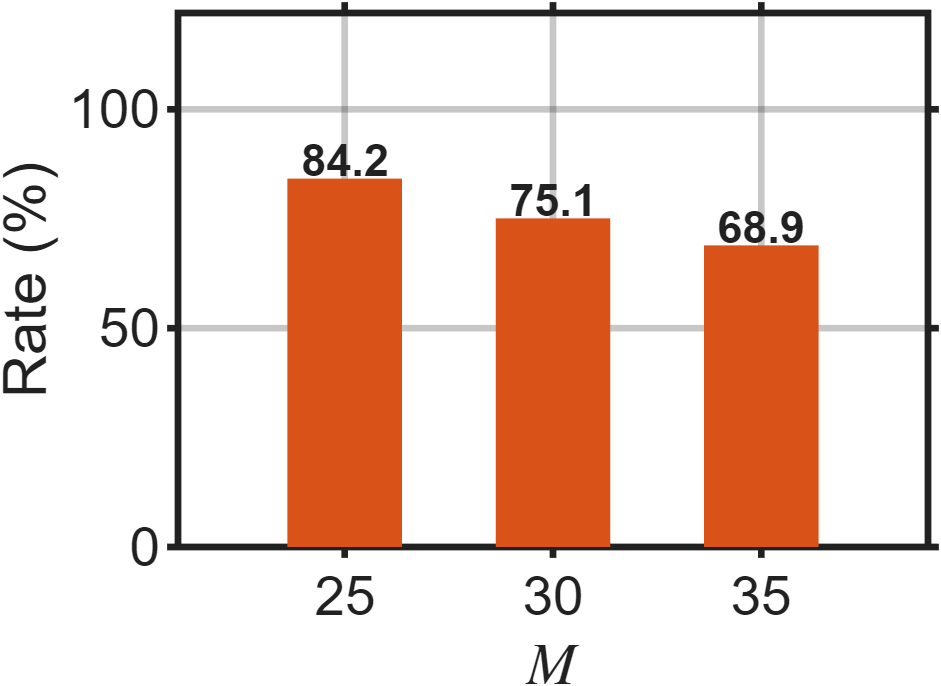}
        \caption{Interception rate under overload ($N=20$)}
        \label{fig:sim2_resilience}
    \end{subfigure}

    \caption{Statistical evaluation of scalability and resilience. (a) Computation time. (b) Interception rates under saturation attacks.}
    \label{fig:exp2_stats}
\end{figure}

The scalability of the algorithm is evaluated by increasing the number of defenders $N$ from 10 to 50. As illustrated in Fig.~\ref{fig:sim2_scalability}, the average local computation time per agent remains nearly constant as the system scale increases, staying between $0.3~\text{ms}$ and $0.5~\text{ms}$. This characteristic indicates the feasibility of the distributed framework for large-scale multi-agent systems. 

Furthermore, the robustness of the partitioning mechanism is evaluated through Monte Carlo simulations under high-load conditions. For each attacker load $M\in\{25,30,35\}$, 50 independent trials are conducted with $N=20$ defenders. In each trial, the defenders are uniformly initialized at a radius of $(R_{in}+R_{out})/2$, while all attackers are simultaneously initialized at random phases and radial distances within $[R_{detect},R_{detect}+5~\text{m}]$. Each trial lasts at most $20~\text{s}$, and the interception rate is calculated as the ratio of the total number of captured attackers to the total number of attackers over all trials. As shown in Fig.~\ref{fig:sim2_resilience}, the system maintains an interception rate of $84.2\%$ at an attacker-to-defender ratio of 1.25 and $68.9\%$ at a ratio of 1.75, demonstrating that the proposed task allocation mechanism remains effective under overloaded conditions.

\subsection{Comparison with a Centralized Assignment Strategy}

This subsection presents a comparative analysis between the proposed algorithm and a centralized strategy that utilizes the Hungarian algorithm for target assignment. The comparison focuses on the operational differences between the proactive positioning logic and the reactive assignment logic.
\begin{figure}[h]
    \centering
    \includegraphics[width=0.92\linewidth]{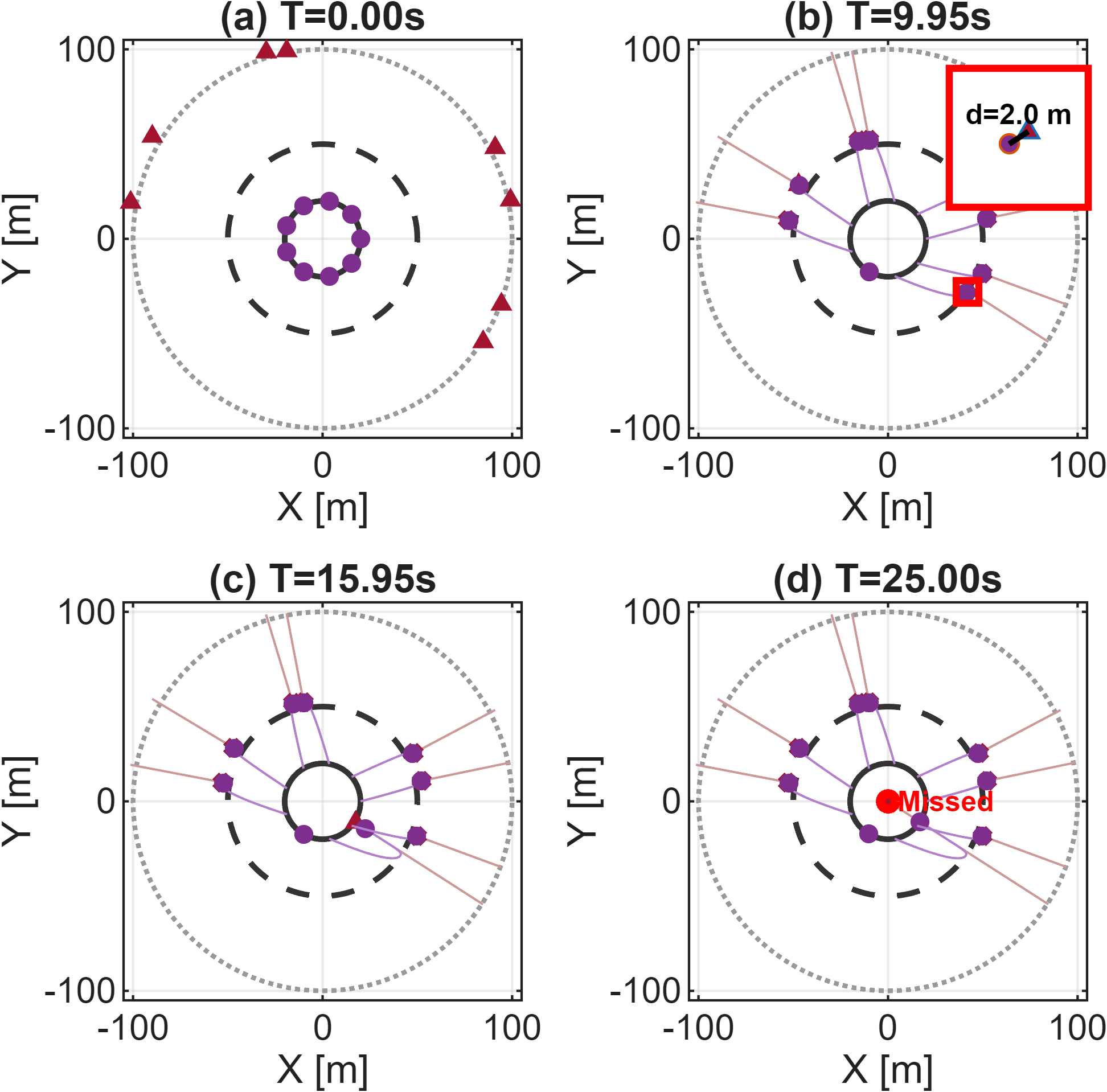}
    \caption{Snapshots of Centralized Strategy. (a) Initialization at protection boundary; (b) Distance gap at closest approach; (c) Tail-chase scenario; (d) Target penetration by attackers.}
    \label{fig:sim3_baseline}
\end{figure}

In the centralized strategy, defenders remain at the protection boundary $R_{in}$ and initiate interception only after receiving assignments upon target detection. Fig.~\ref{fig:sim3_baseline} illustrates the results of this reactive defense. At the moment of closest approach between the defenders and their assigned targets ($t = 9.95~\text{s}$), the minimum separation remains at $2.0~\text{m}$, which exceeds the capture threshold. Due to the lack of prior positioning, defenders enter a tail-chase pursuit but do not close the relative distance before the targets penetrate the defense zone. This outcome highlights the importance of the proactive deployment mechanism for achieving successful interception under kinematic constraints.

\begin{figure}[h]
    \centering
    \begin{subfigure}[b]{0.47\columnwidth}
        \centering
        \includegraphics[width=\textwidth]{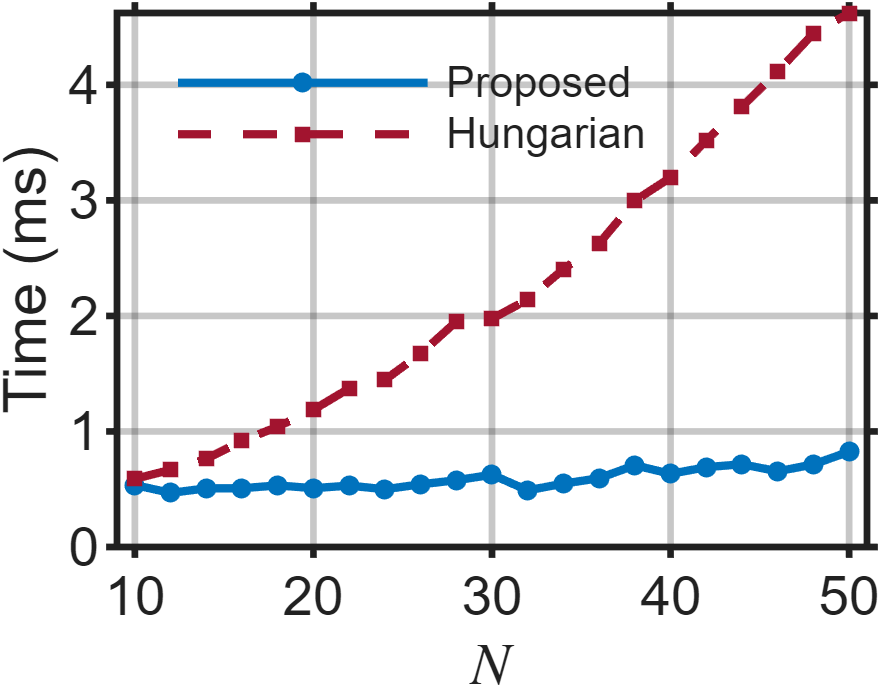}
        \caption{Computational Scalability Comparison}
        \label{fig:sim4_scale}
    \end{subfigure}
    \hfill
    \begin{subfigure}[b]{0.47\columnwidth}
        \centering
        \includegraphics[width=\textwidth]{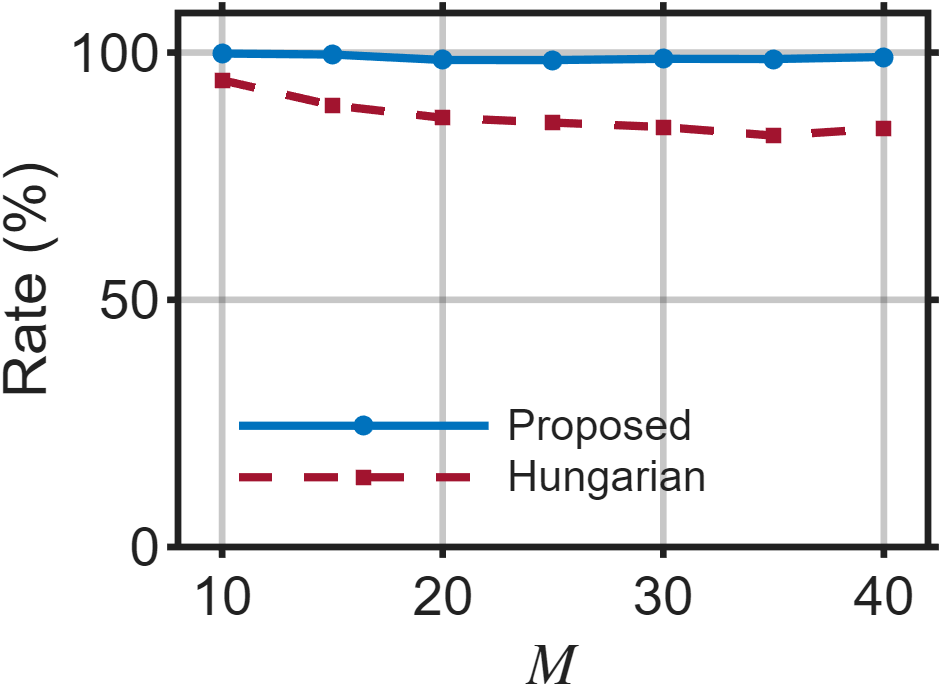}
        \caption{Continuous-arrival interception rates ($N=20$)}
        \label{fig:sim4_resil}
    \end{subfigure}
   \caption{Comparative statistical results. (a) Computation-time comparison. (b) Interception-rate comparison under continuous attacker arrivals.}
    \label{fig:sim4_stats}
\end{figure}

It should be noted that the simulation setting in Fig.~\ref{fig:sim4_stats}(b) differs from the simultaneous saturation-attack scenario in Fig.~\ref{fig:sim2_resilience}. In Fig.~\ref{fig:sim4_stats}(b), attackers are initially separated by $20~\text{m}$ in the radial direction and enter the detection region sequentially. Therefore, this experiment evaluates the ability of the defender network to repeatedly reassign released defenders under continuous arrivals, rather than its instantaneous interception capability under simultaneous overload.

\begin{figure}[!t]
    \centering
    \begin{subfigure}[b]{0.90\columnwidth}
        \centering
        \includegraphics[width=\textwidth]{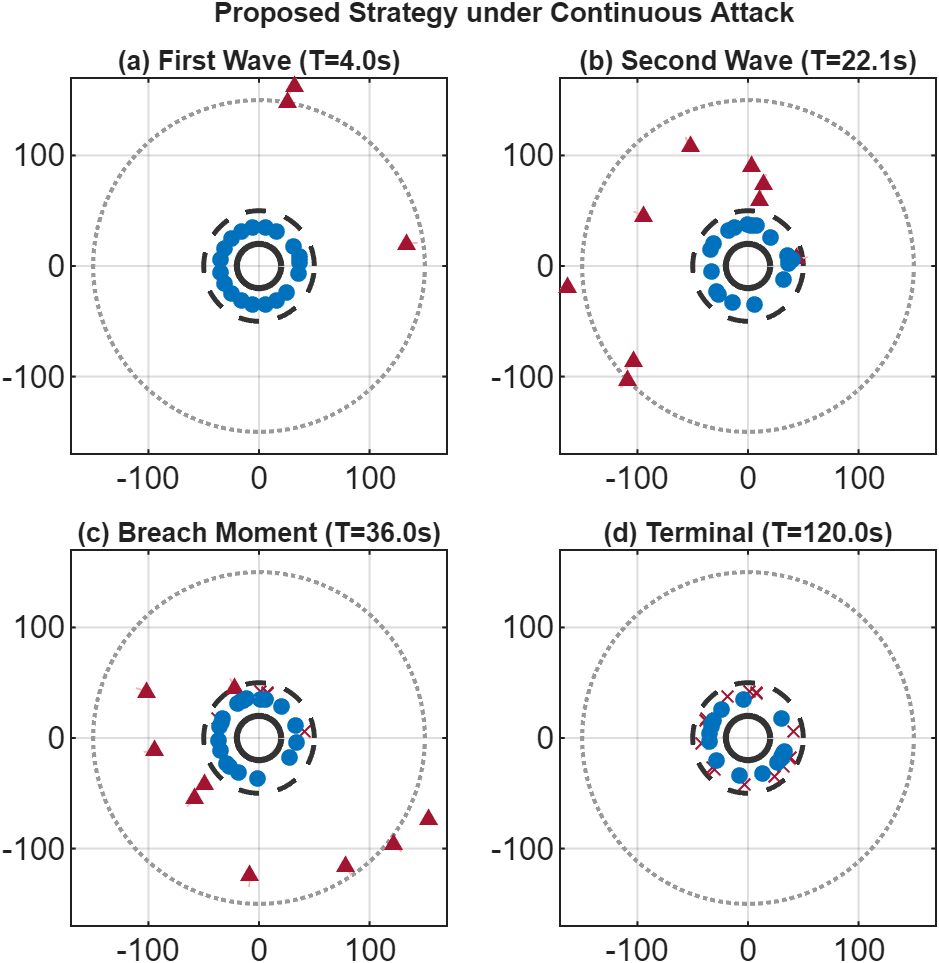}
        \caption{Proposed Method: Defense at $R_{out}$.}
        \label{fig:sim4_vis_prop}
    \end{subfigure}
    \vspace{0.35em}
    \begin{subfigure}[b]{0.90\columnwidth}
        \centering
        \includegraphics[width=\textwidth]{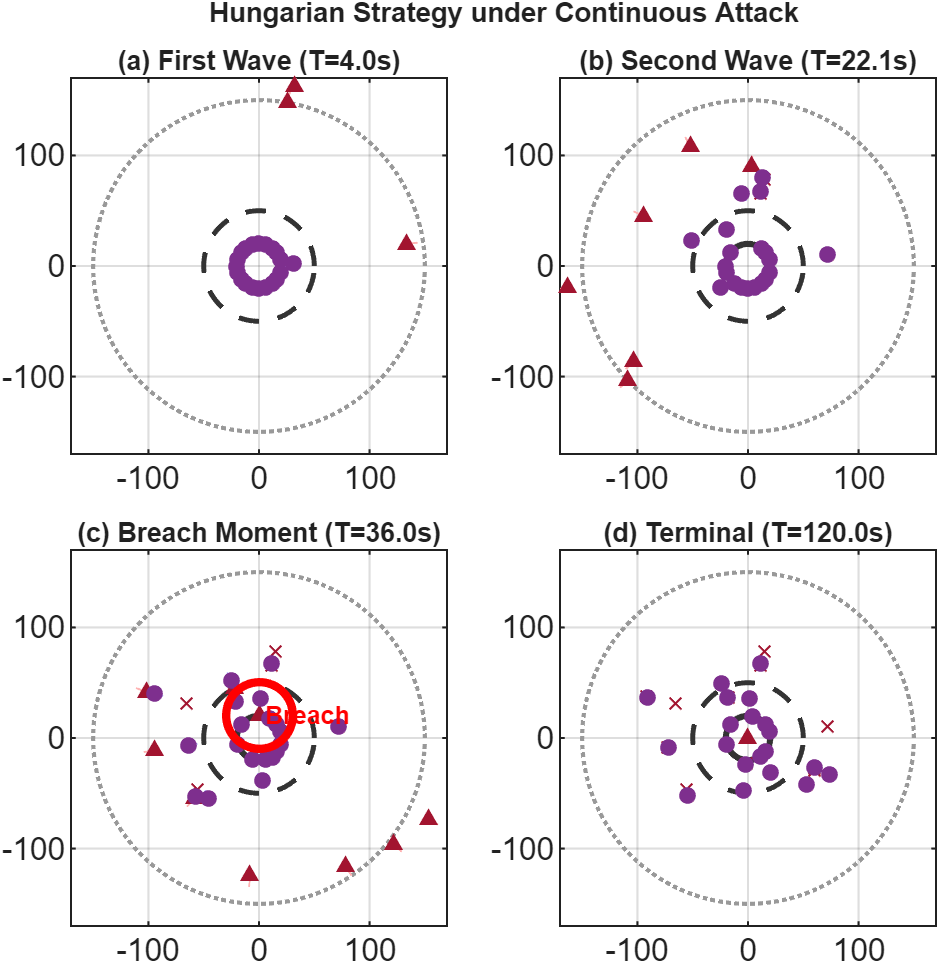}
        \caption{Centralized Strategy: Reactive defense from $R_{in}$.}
        \label{fig:sim4_vis_base}
    \end{subfigure}
    \caption{Visual comparison under continuous attacks ($N=20, M=15$). The proposed method maintains an effective defensive perimeter, while the centralized strategy results in breaches due to reactive positioning.}
    \label{fig:sim4_vis}
\end{figure}

Quantitative results from Monte Carlo simulations further differentiate the two strategies. As depicted in Fig.~\ref{fig:sim4_stats}(a), the proposed distributed algorithm maintains a steady computational load, whereas the centralized strategy exhibits super-linear growth in computation time as the system scales. The interception rates in Fig.~\ref{fig:sim4_stats}(b) indicate that the proposed method experiences a more gradual performance decline as the attacker load increases compared to the centralized strategy. These observations are further supported by the continuous attack scenario shown in Fig.~\ref{fig:sim4_vis}. In this case, the proactive positioning at $R_{out}$ enables the defender network to conduct sequential interceptions. In contrast, the reactive logic of the centralized strategy results in limited effectiveness against a continuous stream of threats.

	\section{Conclusion}\label{sec:4.4}
	
	This paper presented an equitable region partitioning algorithm to address the attacker allocation problem in defender networks. Through mathematical derivation and stability analysis, we demonstrated that, under the specified conditions, the algorithm ensures that each defender is assigned at most one attacker within the annular sectors. For localized one-on-one engagements, we formulated a TAD zero-sum differential game and derived Nash equilibrium strategies solving a Riccati differential equation. Numerical simulations validated the effectiveness and feasibility of the proposed framework, confirming its superior performance in 2D interception scenarios. Future work will explore interception strategies for mobile targets and consider the impact of environmental obstacles.
	
\section*{Acknowledgment}

	The authors used ChatGPT (OpenAI) to polish the language of this manuscript. All authors
	have reviewed the manuscript and approved its final version, and take full responsibility
	for the content of the publication.

	\bibliographystyle{IEEEtran}
	\bibliography{references}

\end{document}